\documentclass[11pt]{article}
\usepackage{amssymb}
\usepackage{amsmath}
\usepackage{amsthm}
\usepackage{latexsym}
\usepackage{amsfonts}
\usepackage{graphicx}
\usepackage{graphics}

\newcommand{\dis}{\displaystyle}
\theoremstyle{plain}
\newtheorem{thm}{Theorem}[]   
\newtheorem{prop}[thm]{Proposition}

\theoremstyle{definition}
\newtheorem{rem}[thm]{Remark}
\newtheorem{exm}[thm]{Example}
\newtheorem*{Proof}{Proof}

\newcommand{\ra}{\;\rightarrow\;}

\newcommand{\al}{\alpha}
\newcommand{\bi}{\beta}
\newcommand{\ga}{\gamma }

\newcommand{\de}{\delta }
\newcommand{\OO} {{\varOmega}}

\newcommand{\e}{\varepsilon }
\newcommand{\f}{\varphi}
\newcommand{\Fi}{\varPhi}

\newcommand{\zi}{\zeta }
\newcommand{\thi}{\theta }

\newcommand{\mi}{\mu }

\newcommand{\ti}{\tau }

\newcommand{\C}{\mathbb{C}}

\newcommand{\ssum}{\sum\limits}

\newcommand{\Ind}{\mbox{Ind}}

\newcommand{\hO}{\widehat{\OO}}

\newcommand{\cu}{{\cal{U}}}

\newcommand{\ld}{\ldots}

\newcommand{\sm}{\smallsetminus}

\newcommand{\qb}{$\quad\blacksquare$}

\begin{document}
\title{One valued primitives and the F. and M. Riesz theorem}
\author{V. Nestoridis}
\date{}
\maketitle\vspace*{-1cm}
\begin{center}
Dedicated to Professor Alain Bernard
\end{center}
\bigskip
\begin{abstract}
For a non-simply connected domain $\OO$ in $\C$ and $f$ a holomorphic function on $\OO$ we prove that $f$ admits one-valued primitives of any order in $\OO$, if and only if, it extends holomorphically in the simply connected envelope of $\OO$. This leads to a generalization of the F. and M. Riesz theorem.
\end{abstract}

{\em AMS classification number}\,: 30E20, 30H10, 30H50, 46J10. \vspace*{0.2cm} 
\section{Introduction}\label{sec1}
\noindent

If $\OO$ is a simply connected domain in $\C$ and $f$ is a holomorphic function in $\OO$, then it is well-known ([1]) that there exists a primitive $F$ of $f$, which is a one-valued holomorphic function satisfying $F'=f$ on $\OO$. To construct $F$ is suffices to fix a point $z_0$ in $\OO$ and for any $z$ in $\OO$ define $F(z)=\dis\int_{\ga_z}f(\zi)d\zi$ for any rectifiable curve $\ga_z$ in $\OO$ starting at $z_0$ and ending at $z$, since by the Cauchy theorem the last integral is independent of the path $\ga_z$ in $\OO$ with the above properties. If $\OO$ is not simply connected the above integrals define in general a multivalued primitive of $f$. However, it may happen that still $F$ is a one-valued holomorphic function in the non-simply connected domain $\OO$. Even more $F$ can still have a one-valued holomorphic primitive on $\OO$ and so on; that is, it may happen that $f$ admits one-valued primitives of any order in the domain $\OO$; equivalently, for every $n=1,2,\ld$ there exists a holomorphic function $\f_n$ on $\OO$, such that the $n$-th derivative $\f^{(n)}_n$ equals $f$ on $\OO$.

In Section \ref{sec1} we prove that the above happens, if and only if, the function $f$ has a holomorphic extension on $\hO$, the least simply connected domain in $\C$ containing $\OO$. Let $X$ be the component of $\infty$ in $(\C\cup\{\infty\})\sm\OO$. Then $\hO=X^c$ and the above condition is that on this simply connected domain the function $f$ has a holomorphic extension. It is also proven in Section \ref{sec2} that the above fact is equivalent to $\dis\int_\ga P(z)f(z)dz=0$ for all analytic polynomials $P$ and all closed polygonal curves $\ga$ in $\OO$.

In Section \ref{sec3} we consider this last condition for all polynomials $P$ on a fixed rectifiable Jordan curve $\ga$ only, where the function $f$ is defined and continuous only on the image $\ga\ast$ of $\ga$. This is equivalent to the fact that $f$ admits one-valued primitives of any order on $\ga$, in some sense to be explained in Section \ref{sec3}. We prove that $f$ has an analytic extension in some sense in the Jordan region $U$ bounded by $\ga$, but we cannot prove that this extension belongs to $A(U)$ and extends continuously on $\overline{U}$. We give an example indicating that this cannot be in general true. The proof uses the classical F. and M. Riesz theorem and we are led to a generalization of it valid for rectifiable Jordan curves. In fact the hypothesis of the F. and M. Riesz is equivalent to the existence of one-valued primitives of any order of a function on a curve.
\section{Main result}\label{sec2}
\noindent

Let $f$ be a holomorphic function on the domain $\OO\subset\C$, such that its primitive is a one-valued holomorphic function $F$. Let $\ga$ be a closed polygonal line in $\OO$. Then $\dis\int_\ga1\cdot f(z)dz=F(\ga^+)-F(\ga^-)=0$, because the endpoints $\ga^+$ and $\ga^-$ of $\ga$ coincide. Thus, the constant function 1 is annilated by the functional defined by the measure $f(z)dz$ on $\ga$. Assume in addition that the holomorphic function $F$ has a one-valued primitive $G$ on $\OO$. Then $\dis\int_\ga zf(z)dz=\dis\int_\ga zdF=\ga^+\cdot F(\ga^+)-\ga^-F(\ga^-)-\dis\int_\ga F(z)dz=G(\ga^+)-G(\ga^-)=0$. Therefore $z^0$ and $z^1$ are annilated by the measure $f(z)dz$ on $\ga$. Continuing in this way by induction we can prove the following.
\begin{prop}\label{prop1}
Let $\OO$ be a domain in $\C$ and $\ga$ a closed polygonal curve in $\OO$. Let $f$ be a holomorphic function on $\OO$, such that, every ordered antiderivative of $f$ is a one-valued holomorphic function on $\OO$; equivalently, we assume that for every $n=1,2,\ld$ there exists a holomorphic function $f_n$ such that $f^{(n)}_n=f$ on $\OO$.

Then every analytic polynomial $P$ is annilated the measure $f(z)dz$ on $\ga$; that is $\dis\int_\ga P(z)f(z)dz=0$.
\end{prop}
\begin{prop}\label{prop2}
Let $V\subset\C$ be a finitely connected domain. Let $V_0,V_1,\ld,V_{k-1}$ be the components of $(\C\cup\{\infty\})\sm V$, $\infty\in V_0$. Let $f$ be a holomorphic function on $V$. For any closed polynomial line $\ga$ in $V$ and every analytic polynomial $P$ we assume that $\dis\int_\ga P(z)f(z)dz=0$. Then, $f$ has a holomorphic extension in the simply connected hull $\widehat{V}$ of $V$, that is on $\widehat{V}=V\cup V_1\cup \cdots \cup V_{k-1}$.
\end{prop}
\begin{Proof}
By Laurent decomposition \cite{C-N-P} $f=f_0+f_1+\cdots+f_{k-1}$, where $f_j\in H(V^c_j)$, for $j=0,1,\ld,k-1$ and for $j=1,\ld,k-1$ we have $\dis\lim_{z\ra\infty}f_j(z)=0$. It suffices to show that $f_1=f_2=\cdots=f_{k-1}=0$. We will do it only for the function $f_1$. This function has a Laurent expansion $f_1(z)=\ssum^\infty_{n=1}\dfrac{a_{-n}}{(n-c)^n}$ for $|z|>R$ where $0<R<+\infty$ is big enough and $c\in V_1$ is fixed.

Let $\ga$ be a closed polygonal curve in $\OO$ surrounding $V_1$ and not $V_2,\ld,V_{k-1}$. That is, $\Ind(\ga,c)=1$ while $\Ind(\ga,c_j)=0$ for $j=2,\ld,k-1$ and $c_j\in V_j$ are fixed \cite{A}. Then $\dis\int_\ga f_j(z)(z-c_j)^ndz=0$ for every $n=0,1,2,\ld$ and $j=2,\ld,k-1$ by Cauchy's theorem. It follows that
\begin{align*}
a_{-n}=&\frac{1}{2\pi i}\int_\ga f_1(z)(z-c)^ndz=\frac{1}{2\pi i}\int_\ga(f_0(z)+f_1(z) \\
&+\cdots+f_{k-1}(z))(z-c)^ndz=\frac{1}{2\pi i}\int_\ga(z-c)^nf(z)dz=0,
\end{align*}
because the polynomial $(z-c)^n$ is annilated by the measure $f(z)dz$ on $\ga$ by assumption. Thus, $f=f_0$ and the proof is completed. \qb
\end{Proof}
\begin{thm}\label{thm3}
Let $\OO\subset\C$ be a domain and $f$ a holomorphic function in $\OO$. Then, the following are equivalent.

(1) The function $f$ admits all orders one-valued primitives; that is, there exists $\f_n\in H(\OO)$, $n=1,2,\ld$ satisfying $\f^{(n)}_n=f$ on $\OO$.

(2) On every closed polygonal line $\ga$ in $\OO$ every complex polynomial $P$ is annilated by the measure $f(z)dz$ on $\ga$; that is $\dis\int_\ga P(z)f(z)dz=0$.

(3) The function $f$ has a holomorphic extension on the simply connected envelope $\widehat{\OO}$ in $\C$ of $\OO$; that is, on the intersection of all simply connected domains $G$, $\OO\subset G\subset\C$. Equivalently on $\widehat{\OO}=X^c$ where $X$ is the component of $\infty$ in $(\C\cup\{\infty\})\sm\OO$.
\end{thm}
\begin{Proof}
(1) implies (2) according to Proposition \ref{prop1}. That (3) implies (1) is also immediate, because $f$ will be holomorphic in the simply connected comain $\widehat{\OO}=X^c\supset\OO$.

The main part of the proof is that (2) implies (3). It is well known that the domain $\OO$ has an exhaustion by compact set $K_m$, $m=1,2,\ld$, where $K_m\subset K^0_{m+1}$, $\bigcup\limits_m K_m=\OO$ and each $K_m$ is the closure of a domain $\OO_m$ of finite connectivity bounded by a finite set of disjoint Jordan curves. We consider $\widehat{\OO_m}$ the simply connected envelope of the domain $\OO_m$. Then according to Proposition \ref{prop2} the function $f$ has a holomorphic extension to $\widehat{\OO}_m$.

The same holds on $\widehat{\OO}_{m+1}\supseteq\widehat{\OO}_m$. The two extensions coincide on $\widehat{\OO}_m$ by analytic continuation, because they coinside on the domain $\OO_m\subset\widehat{\OO}_m\subset\widehat{\OO}_{m+1}$. Thus, we have a holomorphic extension of $f$ on $\bigcup\limits_m\widehat{\OO}_m$. It suffices to show that $\bigcup\limits_m\widehat{\OO}_m$ is the intersection of all simply connected domains in $\C$ containing $\OO$ and that this coincides with $X^c$.

Indeed, the family $(\C\cup\{\infty\})\sm\widehat{\OO}_m$ is a decreasing family of connected compact sets all containing $\infty$. Thus, their intersection is a connected compact set containing $\infty$. This intersection is disjoint from $\OO$ and coincides with $(\C\cup\{\infty\})\sm\bigcup\limits^\infty_{m=1}\widehat{\OO}_m$. Thus, the increasing union of domains $\bigcup\limits_m\widehat{\OO}_m$ is a simply connected domain in $\C$ containing $\OO$. It follows that $\bigcup\limits_m\widehat{\OO}_m\supset\widehat{\OO}$. But every simply connected domain in $\C$ containing $\OO$, it must contain $\widehat{\OO}_m$ for all $m=1,2,\ld$\;.
It follows easily that $\widehat{\OO}=\bigcup\limits_m\widehat{\OO}_m\subset X^c$.
On the other hand let $G\subset\C$ be a simply connected domain containing $\OO$. Let $Y$ be the component of $\infty$ in $(\C\cup\{\infty\})\sm G\subset(\C\cup\{\infty\})\sm\OO$. It follows that $Y\subset X$ and therefore $X^c\subset Y^c=G$. Since this holds for the arbitrary simply connected domain $G$, $\OO\subset G\subset\C$, it follows that $X^c\subset\widehat{\OO}=\cap G$. Thus, $X^c=\widehat{\OO}$ and the proof is\linebreak completed. \qb
\end{Proof}
\begin{rem}\label{rem4}
In view of the equivalence (1)$\Leftrightarrow$(2) of Theorem \ref{thm3} one can easily verify the following.

There exists a holomorphic function $\f_n$ in $\OO$, such that $\f^{(n)}_n=f$ on $\OO$ if and only if the following holds.

For every closed polygonal curve $\ga$ in $\OO$ and every analytic polynomial $P$ with deg$P\le n-1$ we have $\dis\int_\ga P(z)f(z)dz=0$.
\end{rem}
\begin{rem}\label{rem5}
All above results have been stated for closed curves $\ga$ in $\OO$ which are polygonal; however they are valid more generally for closed rectifiable curves in $\OO$; that is, with finite length. It also suffices to consider only rectifiable Jordan curves in $\OO$.
\end{rem}
\section{The F. and M. Riesz theorem for rectifiable Jordan curves}\label{sec3}
\noindent

Let $\OO$ be a doubly connected domain in $\C$. Let $E$ be the bounded component of $\OO^c$ and let us fix a point $a\in E$. Let $\ga$ be a rectifiable Jordan curve in $\OO$ such that $\Ind(\ga,a)=-1$. Let $f$ be a holomorphic function in $\OO$. According to Theorem \ref{thm3} the cunction $f$ has any order's one-valued primitive on $\OO$, if and only if, $\dis\int_\ga P(z)f(z)dz=0$ for all analytic polynomial $P$. And then $f$ has a holomorphic extension in $\OO\cup E$.

We can restrict the construction of the primitives of $f$ on $\ga$. Let $\ga:[\al,\bi]\ra\C$, $\al<\bi$, be a parametrization of the curve $\ga$, $\ga(\al)=\ga(\bi)$. Then for every $z=\ga(t)$, $\al\le t<\bi$ the primitive $G$ of $f$ at $z\in\ga^\ast$ is given by $G(z)=\dis\int_{\ga_z}f(\zi)d\zi$ where $\ga_z:[a,t)\ra\C$ is defined by $\ga_z(\xi)=\ga(\xi)$ for all $\xi\in[a,g]$. The condition that the primitive $G$ is one-valued on $\OO$ is equivalent to the fact that $\dis\int_\ga f(\zi)d\zi=0$ which is the same as $\dis\lim_{t\ra\bi^-}\dis\int_{\ga_{\ga(t)}}f(\zi)d\zi=0$; this is equivalent to the fact that $\dis\int_\ga P(z)f(z)dz=0$ for all analytic polynomials $P$ with deg$P\le0$.

Provided that the first primitive $G$ of $f$ is one valued, in order to construct the second primitive we repeat a similar procedure. Now the fact that the second primitive is one-valued is equivalent to $\dis\lim_{t\ra\bi^-}\dis\int_{\ga_{\ga(t)}}G(\zi)d\zi=0$; and this is equivalent to the fact that $\dis\int_\ga P(z)f(z)dz=0$ for all analytic polynomials $P$ with deg$P\le1$.

Iterating this procedure we conclude that $f$ admits one-valued primitive in $\OO$ of any order, if and only if, $\dis\int_\ga P(z)f(z)dz=0$ for all analytic polynomials $P$. We see that it suffices to examine the situation only on the curve $\ga$ and we could speak about primitives on $\ga$ of functions defined only on $\ga$.

Let us return to the case where $f$ is holomorphic on the doubly connected domain $\OO$. Assume that $E$ the bounded component of $\OO^c$ is a closed Jordan region bounded by a rectifiable Jordan curve $\de$ and that $f$ extends continuously on $\OO\cup\de^\ast$. Then in the previous results we can replace $\ga$ by $\de$. In order to see this, let $\OO_1=(\C\cup\{\infty\})\sm E$ and $\OO_2=\Big\{\dfrac{1}{w-a}:w\in\OO_1\Big\}$, where $a$ is in $E^0$.

Then $\OO_2$ is a Jordan domain bounded by a rectifiable Jordan curve. Let $\Fi:D\ra\OO_2$ be a Riemann map from the open unit disc $D$ onto $\OO_2$. Then $\Fi'$ is in the Hardy Space $H^1$ \cite{D}; it follows that for the curves $\ga_r(t)=\dfrac{1}{\Fi(re^{it})}+a$, $0\le t\le2\pi$, $0<r<1$ and any analytic polynomial $P$ we have
\[
\lim_{r\ra 1^-}\int_{\ga_r}P(z)f(z)dz=\int_{-\de}P(z)f(z)dz.
\]

Therefore, we see that if $\de$ is a rectifiable Jordan curve with a one-sided collar $\OO$ and $f$ is holomorphic on $\OO$ with continuous extension on $\OO\cup\de^\ast$, then the following holds.

The function $f$ admits all orders one-valued primitives on $\OO$, if and only if, it extends holomorphically on $\OO\cup E$, where $E$ is the closed Jordan region bounded by $\de$. And this happens, if and only if, $\dis\int_\de P(z)f(z)dz=0$ for all analytic polynomials $P$. This suggests the following.
\begin{prop}\label{prop6}
Let $\de$ be a rectifiable Jordan curve and let $U$ be the open region bounded by $\de^\ast$. Let $f:\de^\ast\ra\C$ be a continuous function. We assume that $\dis\int_\de P(z)f(z)dz=0$ for all analytic polynomials $P$. (This could be equivalently stated as that $f$ has all orders ``primitives on $\de$'' which are one valued on $\de^\ast$). Then there exists a holomorphic function $\f$ on $U$, such that, for almost all points $j$ in $\de^\ast$ with respect to the arc length measure $|dz|$, the non-tangential limit $\dis\lim_{nt,z\ra \zi\atop z\in U}\f(z)$ exists and equals $f(\zi)$.
\end{prop}
\begin{Proof}
According to Theorem 10.4 of \cite{D}, the function $f$ has an extension $\f$ in $U$ belonging to the class $E^1$ defined in \cite{D}, which satisfies the requirements of Proposition \ref{prop6}.

For completeness we present now another proof.

Let $\Fi:D\ra\cu$ be a Riemann map from the open unit disc $D$ onto $\cu$. Then $\Fi'$ belongs to the Hardy class $H^1$ \cite{D}. We consider on $\de^\ast$ the measure $f(z)dz=d\mi(z)$, which annilates all polynomials $P$. By Mergelyan's theorem \cite{R}, the measure $d\mi$ annilates all elements of the algebra $A(\cu)=\{g:\overline{\cu}\ra\C$ continuous on $\overline{\cu}$, holomorphic in $u=\overline{\cu}^0\}$. We compose with the Riemann map $\Fi:D\ra\cu$, which by the Oswood-Caratheodory theorem \cite{K} extends to a homeomorphism $\Fi:\overline{D}\ra\overline{\cu}$. Then the image of $d\mi$ is a measure $dv$ on the unit circle $T=\partial D$, which annilates all elements of $A(D)$. By the F. and M. Riesz theorem \cite{H} the measure $dv$ is absolutely continuous with respect to Lebesgue measure and $dv=h(e^{i\thi})de^{i\thi}$ on $T$ with $h$ an $H^1$ function on $D$.

Since $d\mi(z)=f(z)dz$ on $\de$, its image composing with $\Fi$ is $f(\Fi(e^{i\thi}))$ multiplied with the image of $dz$ on $\de$ after compositgion with $\Fi$, which is $\Fi'(e^{i\thi})de^{i\thi}+d\ti$, where $d\ti$ is a singular measure on $T$. Thus, the image of $d\mi(z)$ after composition with $\Fi$ is $f(\Fi(e^{i\thi}))\Fi'(e^{i\thi})de^{i\thi}+f(\Fi(e^{i\thi}))dr(e^{i\thi})$.

By the F. and M. Riesz theorem \cite{H} it follows $f(\Fi(e^{i\thi}))d\ti(e^{i\thi})=0$ and $f\circ\Fi\cdot\Fi'=h$ with $h$ in $H^1$ of $D$. Thus, $f\circ\Fi=\dfrac{h}{\phi'}$ almost everywhere on $T=\thi D$ with respect to the Lebesgue measure on $T$. [In fact, since the measure $1\cdot d\ti$ annilates all polynomials, we can repeat everything with $f\equiv 1$. This gives $d\ti\equiv0$].

Since $\Fi'$ is in $H^1$ of $D$ and $\Fi'\not\equiv0$, Jensens formula \cite{H} implies that the non-tangential limits exist almost everywhere and they are non-zero almost everywhere on $T$. Similarily, since $h\in H^1$ of $D$, the non-tangential limits of $h$ exist almost everywhere on $T=\partial D$ with respect to Lebesgue measure. Since $\Fi$ is a univalent holomorphic function on $D$, it follows that $\Fi'(z)\neq0$ for all $z\in D$. Thus, $\dfrac{h}{\Fi'}$ is holomorphic in $D$ and its non-tangential limits on $\thi D$ exist almost everywhere on $T=\thi D$ and are equal to $f\circ\phi$. Then the function $\f=\dfrac{h}{\Fi'}\circ(\Fi^{-1})$ has the desired properties and the proof is\linebreak complete. \qb
\end{Proof}
\begin{exm}\label{exm7}
Since the function $f$ is continuous on $\de^\ast$, naturally comes the question whether the function $\dfrac{h}{\Fi'}\circ(\Fi^{-1})$ belongs to $A(\cu)$; which is equivalent to aks whether $\dfrac{h}{\Fi'}\in A(D)$. We do not have a complete answer to this question. If the Jordan curve $\de$ is such that $|\Fi'(z)|>\e$ for all $z\in D$ and some $\e>0$ independent of $z$, then the answer is affirmative and $\dfrac{h}{\Fi'}\in A(D)$. An example hinting to the negative direction is given by the function $B(z)=\dfrac{z-1}{\exp\dfrac{z+1}{z-1}}$.

The non-tangential limits of it exist on $T-\{1\}$ and define a continuous function on $T$ vanishing it $\zi=1$. However, $B$ is not bounded on $D$ and $B\notin A(D)$.

If $B|T$ had an extension in $A(D)$, this should coincide with $B$, by the reflextion principle or even by Privalor's theorem \cite{K}.
\end{exm}

The weak point of the above example is that the singular inner function $\exp\dfrac{z+1}{z-1}$ does not seem to have a univalent primitive on $D$.
\begin{rem}\label{rem8}
Proposition \ref{prop6} may easily take the following more general form.

Let $\de$ be a rectifiable Jordan curve and $\cu$ the Jordan region bounded by $\de$. Let $d\mi$ be a Borel measure on $\de^\ast$ such that $\dis\int_{\de^\ast}P(z)d\mi(z)=0$ for all analytic polynomials $P$. Then there exists a holomlorphic function $\f$ in $U$ admitting non-tangential limits almost everywhere on $\de^\ast$ with respect to the arc length measure on $\de$ and $d\mi=\f(z)dz$ on $\de$. The function $\f$ belongs to the class $E^1(\cu)$ \cite{D} and conversely for every $\f$ in $E^1(\cu)$ the measure $d\mi=\f(z)dz$ on $\de$ annilates all polynomials; that is $\dis\int_{\de^\ast}P(z)d\mi(z)=0$ for all analytic polynomials. We recall that an equivalent definition of the class $E^1(\cu)$ is the following \cite{D}.

$\f\in E^1(\cu)\Leftrightarrow\f\circ\Fi\cdot\Fi'\in H^1(D)$, where $D$ denotes the open unit disc.
\end{rem}
\begin{thm}\label{thm9}
Let $\cu$ be a Jordan domain bounded by a rectifiable Jordan curve $\de$ and let $\Fi:D\ra\cu$ be a Riemann map from the open unit disc $D$ onto $\cu$. Let $\mi$ be a Borel measure on $\de^\ast$. Then (4) and (5) below are equivalent.

(4) For every analytic polynomial $P$ we have $\dis\int_\de P(z)d\mi(z)=0$

(5) $d\mi(z)=\f(z)dz$, with $\f\in E^1(U)$.
\end{thm}
\begin{Proof}
For the implication (5)$\Rightarrow$(4) we refer to \cite{D}, because the product of the function $\f$ with any polynomial is again in the class $E^1(\cu)$. For the implication (4)$\Rightarrow$(5) we observe that $\dis\int_\de f(z)d\mi(z)=0$ for every function $f:\overline{\cu}\ra\C$ continuous on $\overline{\cu}$ and holomorphic in 4 (that is, $f\in A(\cu)$), because polynomials are dense in $A(\cu)$ by Mergelyan's theorem \cite{R}. Composing with $\Fi$ the image measure $dv$ on $T=\thi D$ of the measure $d\mi$ annilates every function in $A(D)$; in particular every monomial $z^n$, $n=0,1,2,\ld$\;. According to the F. and M. Riesz theorem $dv=h(e^{i\thi})d(e^{i\thi})$ with $h\in H^1(D)$. It follows that $d\mi(z)=h\circ\Fi^{-1}(z)\dfrac{1}{\Fi'(\Fi^{-1}(z))}dz$.

Setting $\f=\dfrac{h\circ\Fi^{-1}}{\Fi'\circ\Fi^{-1}}$ we have $d\mi=\f\cdot dz$ and $\f\in E^1(\cu)$ because $\f\circ\Fi\cdot\Fi'=h$ belongs in $H'(D)$. This completes the proof.  \qb
\end{Proof}
\begin{rem}\label{rem10}
We notice that the measure $dz$ on $\de^\ast$ is the image by composition by $\Fi$ of the measure $\Fi'(e^{i\thi})de^{i\thi}$ on $T$. Since $\de$ is rectifiable we have $\Fi'\in H^1(D)$ \cite{D}. The arc length measure on $\de$ is $|dz|$ and is the image by composition by $\Fi$ of the measure $|\Fi'(e^{i\thi})|d\thi$ on $T$. If $E\subset T$ is a measurable set with zero length then the length of $\Fi(E)$ is $\dis\int_E|\Fi'(e^{i\thi})|d\thi=0$. If $E\subset T$ is a measurable set with strictly positive measure, then the length of $\Fi(E)$ is $\dis\int_E|\Fi'(e^{i\thi})|d\thi>0$; because, if $\dis\int_E|\Fi'(e^{i\thi})|d\thi=0$ it implies $\Fi'\equiv0$ almost everywhere on $E$; this implies $\dis\int_0^{2\pi}\log|\Fi'(e^{i\thi})|d\thi=-\infty$, which contradicts Jensen's inequality \cite{H}, as $\Fi'\in H^1(D)$, $\Fi'\not\equiv0$.
\end{rem}

The above facts can be used to check details in the proof of Theorem \ref{thm9}.
\begin{rem}\label{rem11}
Under the assumptions of Theorem \ref{thm9} we notice that (4) and (5) are equivalent to (6) below.
\end{rem}

Assume that (5) holds. Since $\dis\int_\de\f(z)dz=0$ \cite{D} and $d\mi(z)=\f(z)dz$, it follows that $\mi:\de^\ast\ra\C$ is a continuous one-valued function on $\de^\ast$ and the function $\mi$ takes equal values at the end-points of $\de^\ast$, which coincide. In fact $\mi$ is the restriction on the boundary of the complex primitive $G$ of $\f$ in $\cu$, which belongs to $A(U)$, by Hardy's inequality \cite{H} applied to $h\in H^1(D)$.

Fix $z_0$ in $\de^\ast$. For any $z$ in $\de^\ast$ we consider a curve $\de_z$ in $\de^\ast$ starting at $z_0$ and ending at $z$. In fact $\de_z$ can be seeing as a restriction of the periodic extension of $\de$. We set $F(z)=\dis\int_{\de_z}\mi(\zi)d\zi$. This function is one-valued on $\de^\ast$ and it corresponds to the boundary function of the complex primitive in  $\cu$ of the function $G$, which belongs to $A(\cu)$.

Continuing in this waqy we see that the following hold.

(6) The function $\mi$ has one-valued continuous $dz$-primitives on $\de^\ast$ in the sense of $dz$-integrals of any order $n=0,1,2,3,\ld$\;.

Assume now that (6) holds. We will show that (4) also holds.

Clearly (6) implies $\dis\int_\de1d\mi(z)=0$. Also $F$ is a one-valued function on $\de^\ast$; then $0=\dis\int_\de\mi(\zi)d\zi=\zi\mi(\zi)|_{\de^\ast}-\dis\int_\de \zi d\mi(\zi)=0-\dis\int_\de \zi d\mi(\zi)$; this implies $\dis\int_\de \zi d\mi(\zi)=0$.

Continuing in this way we see that $\dis\int_\de \zi^nd\mi(\zi)=0$ for all $n=0,1,2,\ld$ which implies (4) by the linearity of the integral. We also notice that in the last argument integrations by parts may be replaced by the use of Fubini's theorem.

We also notice that if the rectifiable Jordan curve $\de$ satisfies some assumption and $F$ is a one-valued $dz$-primitive of the continuous function $G$ on $\de^\ast$, then
\[
\lim_{z\ra z_0\atop z\in\de^\ast}\frac{F(z)-F(z_0)}{z-z_0}=G(z_0) \ \ \text{for all} \ \ z_0\in\de^\ast.
\]
Indeed $F(z)-F(z_0)=\dis\int_{d_{z_0,z}}G(j)dj$ where $\de_{z_0,z}$ is the ``shortest'' arc of the curve $\de$ starting at $z_0$ and ending at $z$. It follows that
\[
\frac{F(z)-F(z_0)}{z-z_0}-G(z_0)=\frac{1}{z-z_0}\int_{\de_{z_0,z}}[G(\zi)-G(z_0)]d\zi.
\]
Taking absolute values this is bounded above by $\dfrac{{\text{length of}\,\de_{z_0,z}}}{|z-z_0|}\cdot\dis\sup_{\zi\in\de_{z0,z}}|G(\zi)-G(z_0)|$.

If we assume that there exists a constant $C_{z_0}<+\infty$ such that (length of $\de_{z_0,z}$)$\le C_{z_0}\cdot|z-z_0|$ for all $z\in\de^\ast$ close to $z_0$, then we take $\Big|\dfrac{F(z)-F(z_0)}{z-z_0}-G(z_0)\Big|\le C_{z_0}\sup_{\zi\in\de_{z_0,z_0}}|G(\zi)-G(z_0)|\ra0$, as $z\ra z_0$, $(z\in\de^\ast)$, by the continuity of $G$.

In particular this holds if the rectifiable Jordan curve $\de$ is a chord-arc curve; that is, if there exist a constant $C<+\infty$ so that (length of $\de_{z_0,z}$)$\le C\cdot|z-z_0|$ for $z,z_0$ on $\de^\ast$.\medskip\\
\noindent
{\bf Acknowledgment.} I would like to thank John Pardon for his essential help in the proof of Theorem \ref{thm3} and D. Gatzouras, G. Koumoullis, N. Papadatos, Ch. Panagiotis and T. Hatziafratis for helpful discussions.
V. Nestoridis, \\
National and Kapodistrian University of Athens \\
Department of Mathematics\\
Panepistemiopolis \\
157 84 Athens\\
Greece\\
e-mail:vnestor@math.uoa.gr

\end{document}